\documentclass{article}
\newcommand{\Z}{\ensuremath{\mathbf Z}}
\newcommand{\N}{\ensuremath{ \mathbf N }}
\newtheorem{theorem}{Theorem}
\newcommand{\bt}{\begin{theorem}}
\newcommand{\et}{\end{theorem}}
\newcommand{\pf}{{\bf Proof}.\ }
\newcommand{\bq}{\begin{eqnarray*}}
\newcommand{\eq}{\end{eqnarray*}}
\newcommand{\be}{\begin{eqnarray}}
\newcommand{\ee}{\end{eqnarray}}
\newcommand{\beq}{\begin{equation}}
\newcommand{\eeq}{\end{equation}}
\newcommand{\benum}{\begin{enumerate}}
\newcommand{\eenum}{\end{enumerate}}
\newcommand{\ba}{\begin{array}}
\newcommand{\ea}{\end{array}}
\newcommand{\card}{\mbox{card}}

\begin{document}
\title{Unique representation bases for the integers}
\author{Melvyn B. Nathanson\thanks{This work was supported
in part by grants from the NSA Mathematical Sciences Program
and the PSC-CUNY Research Award Program.}\\
Department of Mathematics\\
Lehman College (CUNY)\\
Bronx, New York 10468\\
Email: nathansn@alpha.lehman.cuny.edu}
\maketitle

\begin{abstract}
Arbitrarily sparse sets $A$ of integers are constructed with the
property that every integer can be represented uniquely
in the form $n = a+a',$ where $a,a' \in A$ and $a \leq a'$.
Some related open problems are stated.
\end{abstract}

\section{Additive bases for the integers}
Let $A$ be a set of integers, and let $r_A(n)$ denote the number
of representations of $n$ in the form $n = a+a',$ where $a,a' \in A$
and $a \leq a'.$  The function $r_A(N)$ is called the {\em representation
function} of the set $A$.  An unsolved problem of Erd\H os and Tur\' an
states that if $A$ is a subset of the semigroup \N$_0$ of nonnegative integers
and $r_A(n) \geq 1$ for all sufficiently large integers $n$,
then the representation function $r_A(n)$ is unbounded.
On the other hand, it is known
(in the oral tradition, but, perhaps, not the written
tradition of additive number theory)
that the group of integers \Z\ contains sets $A$ with the property
that $r_A(n) \geq 1$ for all $n \in \Z$ and $r(n)$ is bounded.

A set $A$ of integers is called an {\em additive basis for the integers}
if $r_A(n) \geq 1$ for all $n \in \Z$, and a {\em unique representation basis}
if $r_A(n) = 1$ for all $n \in \Z$.
The purpose of this paper is to construct a family of arbitrarily
sparse unique representation bases for \Z.
When a greedy algorithm is used in this construction, we obtain
a unique representation basis $A$ whose growth is logarithmic
in the sense that the number of elements $a \in A$ with $|a| \leq x$
is bounded above and below by constant multiples of $\log x$.
In the last section of this paper we state some open problems
suggested by the additive bases that we have constructed.

\section{Bases with arbitrarily slow growth}
For sets $A$ and $B$ of integers and for any integer $c$,
we define the {\em sumset}
\[
A+B = \{a+b: a \in A, b \in B\}
\]
and the {\em translation}
\[
A + c = \{a+c: a \in A\}.
\]
For the sumset
\[
2A = A+A = \{a+a': a,a' \in A\},
\]
we have the {\em representation function}
\[
r_A(n) = \card\{(a,a')\in A\times A : a \leq a' \mbox{ and } a+a' = n\}.
\]
The {\em counting function} for the set $A$ is
\[
A(y,x) = \card\{a\in A: y \leq a \leq x\}.
\]
In particular, $A(-x,x)$ counts the number of integers $a\in A$
such that $|a| \leq x.$

\bt   \label{z:t1}
Let $f(x)$ be a function such that $\lim_{x\rightarrow\infty} f(x) = \infty.$
There exists an additive basis $A$ for the group \Z\ of integers such that
\[
r_A(n) = 1 \qquad\mbox{for all $n \in \Z$,}
\]
and
\[
A(-x,x) \leq f(x)
\]
for all sufficiently large $x$.
\et

\pf
We shall construct an ascending sequence of finite sets
$A_1 \subseteq A_2 \subseteq A_3 \subseteq \cdots$
such that
\[
|A_k| = 2k \qquad\mbox{for all $k \geq 1$},
\]
\[
r_{A_k}(n) \leq 1 \qquad\mbox{for all $n \in \Z,$}
\]
and
\[
r_{A_{2k}}(n) = 1 \qquad\mbox{for all $n$ such that $|n| \leq k$.}
\]
It follows that the infinite set
\[
A = \bigcup_{k=1}^{\infty} A_k
\]
is a unique representation basis for the integers.

We construct the sets $A_k$ by induction.
Let $A_1 = \{0, 1\}$.  We assume that for some $k \geq 1$ we have
constructed sets
\[
A_1 \subseteq A_2 \subseteq \cdots \subseteq A_k
\]
such that $|A_k| = 2k$ and
\[
r_{A_k}(n) \leq 1 \qquad\mbox{for all $n \in \Z$.}
\]
We define the integer
\[
d_k = \max\{|a| : a\in A_k\}.
\]
Then
\[
A_k \subseteq [-d_k,d_k]
\]
and
\[
2A_k \subseteq [-2d_k, 2d_k].
\]
If both numbers $d_k$ and $-d_k$ belong to the set $A_k$,
then, since $0 \in A_1 \subseteq A_k$ and $d_k \geq 1,$, 
we would have the following two representations of 0 in the sumset $2A_k$:
\[
0 = 0 + 0 = (-d_k)+ d_k.
\]
This is impossible, since $r_{A_k}(0) \leq 1$, hence only one of the
two integers $d_k$ and $-d_k$ belongs to the set $A_k$.
It follows that if $d_k \not\in A_k,$ then
\[
\{2d_k,2d_k-1\} \cap 2A_k = \emptyset,
\]
and if $-d_k \not\in A_k,$ then
\[
\{-2d_k,-(2d_k-1)\} \cap 2A_k = \emptyset.
\]
Define the integer $b_k$ by
\[
b_k = \min\{|b| : b\not\in 2A_k\}.
\]
Then
\[
1 \leq b_k \leq 2d_k-1.
\]

To construct the set $A_{k+1}$, we choose an integer $c_k$
such that
\[
c_k \geq d_k.
\]
If $b_k \not\in 2A_k$, let
\[
A_{k+1} = A_k \cup \{b_k + 3c_k, -3c_k\}.
\]
We have
\[
b_k = (b_k+3c_k) + (-3c_k) \in 2A_{k+1}.
\]
If $b_k \in 2A_k$, then $-b_k \not\in 2A_k$ and we let
\[
A_{k+1} = A_k \cup \{-(b_k + 3c_k), 3c_k\}.
\]
Again we have
\[
-b_k = -(b_k+3c_k) + 3c_k \in 2A_{k+1}.
\]
Since
\[
d_k < 3c_k < 3c_k+b_k,
\]
it follows that $|A_{k+1}| = |A_k|+2 = 2(k+1)$.
Moreover,
\[
d_{k+1} = \max\{|a| : a\in A_{k+1}\} = b_k + 3c_k.
\]

For example, since $A_1 = \{0,1\}$ and $2A_1 = \{0,1,2\}$, it follows
that $d_1 = b_1 = 1.$  For $c_1 \geq 1$ we have
\[
A_2 = \{-(1+3c_1),0,1,3c_1\}.
\]
Then
\[
2A_2 = \{-(2+6c_1), -(1+3c_1), -3c_1, -1, 0, 1, 2, 3c_1, 1+3c_1, 6c_1  \}
\]
and  $d_2 = 1+3c_1$ and $b_2 = 2$.

We can assume that $b_k \not\in 2A_k$, hence
$A_{k+1} = A_k \cup \{b_k + 3c_k, -3c_k\}$.
(The argument in the case $b_k \in 2A_k$ and $-b_k \not\in 2A_k$
is similar.)
We shall show that the sumset $2A_{k+1}$ is the disjoint union
of the following four sets:
\[
2A_{k+1} = 2A_k
\cup \left( A_k + b_k + 3c_k \right)
\cup \left( A_k - 3c_k \right)
\cup \{b_k, 2b_k + 6c_k,-6c_k\}.
\]
If $u \in 2A_k$, then
\[
-2c_k \leq -2d_k \leq u \leq 2d_k \leq 2c_k.
\]
Suppose that $v  = a + b_k + 3c_k \in A_k + b_k + 3c_k$,
where $a \in A_k.$
The inequalities
\[
-c_k \leq -d_k \leq a \leq d_k \leq c_k
\]
and
\[
1 \leq b_k \leq 2d_k-1 \leq 2c_k-1
\]
imply that
\[
2c_k + 1 \leq v \leq 6c_k - 1 < 2b_k + 6c_k.
\]
Similarly, if $w = a - 3c_k \in A_k - 3c_k$, then
\[
-6c_k < -4c_k \leq w \leq -2c_k.
\]
These inequalities imply that the sets $2A_k$,
$A_k + b_k + 3c_k$, $A_k - 3c_k$, and $2\{b_k + 3c_k,-3c_k\}$
are pairwise disjoint, unless $c_k = d_k$ and
$-2d_k \in 2A_k \cap 2(A_k - 3d_k)$.
If $-2d_k \in 2A_k$, then $-d_k \in A_k$.
If $-2d_k \in 2(A_k-3d_k)$, then $d_k \in A_k.$
This is impossible, howev
er, because
the set $A_k$ does not contain both integers $d_k$ and $-d_k$.

Since the sets $A_k + b_k + 3c_k$ and $A_k - 3c_k$ are
translations, it follows that
\[
r_{A_{k+1}}(n) \leq 1 \qquad\mbox{for all integers $n$}.
\]

Let $A = \cup_{k=1}^{\infty} A_k$.
For all $k \geq 1$ we have $2 = b_2 \leq b_3 \leq \cdots$
and $b_k < b_{k+2}$, hence $b_{2k} \geq k+1$.
Since $b_{2k}$ is the minimum of the absolute values of the integers
that do not belong to $2A_{2k}$, it follows that
\[
\{-k,-k+1,\ldots,-1,0,1,\ldots, k-1, k\} \subseteq 2A_{2k} \subseteq 2A
\]
for all $k \geq 1$, and so $A$ is an additive basis.
If $r_A(n) \geq 2$ for some $n$,
then $r_{A_k}(n) \geq 2$ for some $k$, which is impossible.
Therefore, $A$ is a unique representation basis for the integers.

In the construction of the set $A_{k+1}$, the only constraint on the
choice of the number $c_k$ was that $c_k \geq d_k$.
Given a function $f(x)$ that tends to infinity,
we use induction to  construct a sequence of integers
$\{c_k\}_{k=1}^{\infty}$
such that $A(-x,x) \leq f(x)$ for all $x \geq c_1.$
We observe that
\bq
A(-x,x) 
& = & A_{k+1}(-x,x) \\
& = & \left\{\ba{ll}
2k & \mbox{for $d_k\leq x < 3c_k$,}  \\
2k+1 & \mbox{for $3c_k\leq x < b+k+3c_k = d_{k+1}$.}
\ea\right.
\eq
We begin by choosing an integer $c_1 \geq d_1$ such that
\[
f(x) \geq 4 \qquad \mbox{for $x \geq c_1$.}
\]
Then
\[
A(-x,x)  \leq 4 \leq f(x) \qquad \mbox{for $c_1 \leq x \leq d_2$.}
\]
Let $k \geq 2$, and suppose we have selected an integer $c_{k-1} \geq d_{k-1}$ such that
\[
f(x) \geq 2k \qquad \mbox{for $x \geq c_{k-1}$}
\]
and
\[
A(-x,x) \leq f(x) \qquad \mbox{for $c_1 \leq x \leq d_k$.}
\]
There exists an integer $c_k \geq d_k$ such that
\[
f(x) \geq 2k+2 \qquad \mbox{for $x \geq c_k$}
\]
Then
\[
A(-x,x) = 2k \leq f(x) \qquad \mbox{for $d_k \leq x < 3c_k$}
\]
and
\[
A(-x,x) \leq 2k+2 \leq f(x) \qquad \mbox{for $3c_k \leq x \leq d_{k+1}$,}
\]
hence
\[
A(-x,x) \leq f(x) \qquad \mbox{for $c_1 \leq x \leq d_{k+1}$.}
\]
It follows that
\[
A(-x,x) \leq f(x) \qquad \mbox{for all $x \geq c_1$.}
\]
This completes the proof.

\section{Bases with logarithmic growth}
In Theorem~\ref{z:t1} we constructed unique representation bases
whose counting functions tend slowly to infinity.
It is natural to ask if there exist unique representation bases
that are dense in the sense that their counting functions tend
rapidly to infinity.  In the following theorem we use the previous
algorithm to construct a unique representation basis $A$ whose
counting function $A(-x,x)$ has order of magnitude $\log x$.

\bt                        \label{z:t2}
There exists a unique representation basis $A$ for the integers such that
\[
\frac{2\log x}{\log 5} + 2\left( 1 - \frac{\log 3}{\log 5} \right)
\leq A(-x,x) \leq \frac{2\log x}{\log 3} + 2.
\]
\et

\pf
We apply the method of Theorem~\ref{z:t1} with 
\[
c_k = d_k \qquad \mbox{for all $k \geq 1$.}
\]
This is essentially a greedy algorithm
construction, since at each iteration we choose the
smallest possible value of $c_k$.
It is instructive to compute the first few sets $A_k$.
Since
\[
A_1 = \{0,1\} \qquad\mbox{and} \qquad 2A_1 = \{0,1,2\},
\]
we have $b_1 = 1$ and $c_1 = d_1 = 1$.  Then
\[
A_2 = \{-4,0,1,3\}
\]
and
\[ 
2A_2 = \{-8,-4,-3,-1, 0,1,2,3,4,6\},
\]
hence $b_2 = 2$, $c_2 = d_2 = 4$.
The next iteration of the algorithm produces the sets
\[
A_3 = \{-14,-4,0,1,3,12\}
\]
and
\bq
2A_3 & = & \{-28,-18,-14,-13,-11,-8,-4,-3\} \\
& & \cup \{-2, -1, 0,1,2,3,4,6,8,12,15,24 \},
\eq
we obtain $b_3 = 5$, $c_3 = d_3 = 28,$ and
\[
A_4 = \{-84,  -14, -4, 0, 1,3,12,89 \}.
\]

We shall compute upper and lower bounds for the counting function
$A(-x,x)$.
For $k \geq 1$ we have $1 \leq b_k \leq 2c_k-1$ and
$c_{k+1} = 3c_k + b_k$, hence
\[
3c_k + 1 \leq c_{k+1} \leq 5c_k -1.
\]
Since $c_1 = 1,$ it follows by induction on $k$ that
\[
\frac{3^k -1}{2} \leq c_k \leq \frac{3\cdot 5^k + 5}{20}
\]
and so
\[
\frac{\log c_k}{\log 5} + 1
\leq \frac{\log\left( (20c_k-5)/3)\right)}{\log 5}
\leq k
\leq \frac{\log(2c_k+1)}{\log 3}
\leq \frac{\log c_k}{\log 3} + 1
\]
for all $k \geq 1$.
We obtain an upper bound for $A(-x,x)$ as follows.
If $c_k \leq x < 3c_k$, then
\bq
A(-x,x) & = & A_k(-x,x)\\
& = & 2k \\
& \leq & \frac{2\log c_k}{\log 3} + 2  \\
& \leq & \frac{2\log x}{\log 3} + 2.
\eq
If $3c_k \leq x < c_{k+1}$, then
\bq
A(-x,x) & = & A_{k+1}(-x,x) \\
& = & 2k + 1\\
& \leq & \frac{2\log c_k}{\log 3} + 3 \\
& \leq & \frac{2\log (x/3)}{\log 3} + 3 \\
& \leq & \frac{2\log x}{\log 3}+1.
\eq
Therefore,
\[
A(-x,x) \leq \frac{2\log x}{\log 3} + 2 \quad\mbox{for all $x \geq 1$.}
\]
We obtain a lower bound for $A(-x,x)$ similarly.
If $c_k \leq x < 3c_k$, then
\bq
A(-x,x) & = & 2k \\
& \geq & \frac{2\log c_k}{\log 5} + 2 \\
& \geq & \frac{2\log (x/3)}{\log 5} + 2 \\
& \geq & \frac{2\log x}{\log 5} + 2\left( 1- \frac{\log 3}{\log 5}\right) \\
& = & \frac{2\log x}{\log 5} + 0.63\ldots.
\eq
If $3c_k \leq x < c_{k+1}$, then, since
\[
c_{k+1} = d_{k+1} = b_k+3c_k \leq 5c_k-1,
\]
we have 
\bq
A(-x,x) & = & 2k + 1\\
& \geq & \frac{2\log c_k}{\log 5} + 3  \\
& > & \frac{2\log (x/5)}{\log 5} + 3  \\
& = & \frac{2\log x}{\log 5} + 1  \\
\eq
Therefore,
\[
A(-x,x) \geq \frac{2\log x}{\log 5} +  2\left( 1- \frac{\log 3}{\log 5}\right)
\quad\mbox{for all $x \geq 1$.}
\]
This completes the proof of the Theorem~\ref{z:t2}.

\section{Heuristics and open problems}
Let $A$ be a set of integers.
If $A$ is a unique representation basis for \Z, or, more generally,
if $A$ is a set of integers with a bounded representation function, 
then $A(-x,x) \ll \sqrt{x}$.
The following simple result gives an explicit upper bound.

\bt                                    \label{z:t3}
Let $A$ be a nonempty set of integers
such that the representation function of $A$ is bounded.
If $r_A(n) \leq r$ for all $n$, then
\[
A(-x,x) \leq \sqrt{8rx}
\]
for all $x \geq r.$
\et

\pf
Let $A(-x,x) = k$.
The number of ordered pairs $(a,a') \in A\times A$
with $-x \leq a \leq a' \leq x$ is exactly $(k^2+k)/2$.
For each of these ordered pairs we have $-2x \leq a+a' \leq 2x$.
For each integer $n \in [-2x,2x]$ there are at most $r$
such pairs $(a,a')$ with $a+a' = n$, and so
\[
\frac{k^2+k}{2} \leq r(4x+1).
\]
It follows that
\[
A(-x,x) = k \leq \sqrt{8rx  + \frac{8r+1}{4}} - \frac{1}{2} \leq \sqrt{8rx}
\]
for $x \geq r$.
This completes the proof.

Theorem~\ref{z:t3} has a natural analogue for sets $A$ nonnegative 
integers.

\bt                                          \label{z:t4}
Let $A$ be a set of nonnegative integers
such that every sufficiently large integer can be represented
as the sum of two elements of $A$.
If $r_A(n) \geq 1$ for all $n > n_0$, then
\[
A(0,x) \geq 2\sqrt{x}-1
\]
for all $x \geq n_0^2.$
If $A$ is a set of nonnegative integers such that
$r_A(n) \leq r$ for all $n \geq 0$, then
\[
A(0,x) \leq 2\sqrt{rx}
\]
for all $x \geq 1.$
\et

\pf
Let $A(0,x) = k$.
Suppose that $r_A(n) \geq 1$ for all $n > n_0$.
The number of ordered pairs $(a,a') \in A\times A$
with $0 \leq a \leq a' \leq x$ is exactly $(k^2+k)/2$.
For each such pair we have $0 \leq a+a' \leq 2x$.
For each integer $n$ with $n_0 < n \leq 2x$ there is at least
one pair $(a,a')$ with $a+a' = n$, and so
\[
\frac{k^2+k}{2} \geq 2x-n_0.
\]
It follows that
\[
A(0,x) = k \geq \sqrt{4x - 2n_0 + \frac{1}{4}} - \frac{1}{2}
\geq 2\sqrt{x}-1
\]
for $x \geq n_0^2$.

Suppose that $r_A(n) \leq r $ for all $n \geq 0$.
If $a,a' \in A$ and $0 \leq a,a' \leq x$, then
$0 \leq a+a' \leq 2x.$
Since $r_A(0) \leq 1$ and $r_A(1) \leq 1$, it follows,
as in the proof of Theorem~\ref{z:t3}, that
\[
\frac{k^2+k}{2} \leq r(2x-1) + 2,
\]
and so
\[
A(0,x) = k \leq \sqrt{4rx  + \frac{17-8r}{4}} - \frac{1}{2} \leq 2\sqrt{rx}
\]
for $x \geq 1$.
This completes the proof.

A set $A$ of nonnegative integers is called a {\em basis} 
(resp. an {\em asymptotic basis}) if every (resp. every sufficiently large) 
integer can be represented as the sum of two elements of $A$.
By Theorem~\ref{z:t1}, there exist arbitrarily sparse sets of integers
that are unique representation bases for \Z.  
On the other hand, by Theorem~\ref{z:t4}, 
a set $A$ of nonnegative integers that is a 
basis or asymptotic basis for the set of nonnegative integers
must have a counting function $A(0,x)$ that grows at least as fast as $\sqrt{x},$
and if the representation function of $A$ is bounded, then $A(0,x)$ cannot grow 
faster than a constant multiple of $\sqrt{x}$.
This phenomenon can be interpreted as follows:  If $0 \leq n \leq N$, 
then there are infinitely many pairs $(a,a')$ of integers
whose sum is $n$, and the summands $a$ and $a'$ can be arbitrarily large in absolute value.
On the other hand, if $a$ and $a'$ are constrained to be nonnegative integers, 
then they must be 
chosen from the finite number of integers in the bounded interval $[0,N]$.
If $A$ is an asymptotic basis, then $A$ is forced 
to contain many numbers in the interval $[0,N]$, 
and this increases the probability that some number has many representations.
This phenomenon may underlie the Erd\H os-Tur\' an conjecture.

Theorem~\ref{z:t1} asserts that a unique representation
basis $A$ for the integers can be arbitrarily sparse, 
while Theorem~\ref{z:t3} states that $A$ cannot be too dense,
since $A(-x,x) \ll \sqrt{x}$.
In Theorem~\ref{z:t2} we constructed a unique representation basis
such that $A(-x,x) \geq (2/\log 5)\log x + 0.63$.
It is not known what functions can be lower bounds for
counting functions of unique representation bases.
Here are some unsolved problems on this theme.

\benum
\item
For each real number $c > 2/\log 5$,
does there exist a unique representation basis $A$
such that $A(-x,x) \geq c\log x$ for all sufficiently large $x$?

\item
Does there exist a unique representation basis $A$ such that
\[
\lim_{x\rightarrow\infty}\frac{A(-x,x)}{\log x} = \infty?
\]

\item
Does there exist a number $\theta > 0$
and a unique representation basis $A$ such that
$A(-x,x) \geq x^{\theta}$ for all sufficiently large $x$?

\item
Does there exist a number $\theta < 1/2$ such that
$A(-x,x) \leq x^{\theta}$ for every unique representation
basis $A$ and for all sufficiently large $x$?
\eenum

\end{document}